\documentclass[12pt]{amsart}
\usepackage{latexsym, amssymb, amsmath, amscd}
 \usepackage{eucal}

    \newtheorem{rema}{Remark}[section]
    \newtheorem{propo}[rema]{Proposition}

   \newtheorem{theo}[rema]{Theorem}
   \newtheorem{def-theo}[rema]{Definition-Theorem}
 \newtheorem{conj}[rema]{Conjecture}
   \newtheorem{defi}[rema]{Definition}
    \newtheorem{lemma}[rema]{Lemma}
    \newtheorem{corol}[rema]{Corollary}
     
  \newtheorem{rmk}[rema]{Remark}

	\newcommand{\p}{\partial}

 \newcommand{\pf}{{\it Proof:}\hspace{2ex}}
 
 \newcommand{\epfv}{\hspace{1em}$\Box$\vspace{1em}}

\newcommand{\bC}{{\mathbb C}}
\newcommand{\bZ}{{\mathbb Z}}
\newcommand{\bR}{{\mathbb R}}
\newcommand{\bQ}{{\mathbb Q}}
\newcommand{\bN}{{\mathbb N}}

\newcommand{\cA}{{\mathcal A}}

\newcommand{\I}{{\operatorname I}}
\newcommand{\poly}{\operatorname{Poly}}
\newcommand{\supp}{\operatorname{Supp}}

\newcommand{\cB}{{\mathcal B}}

\newcommand{\cH}{{\mathcal H}}
\newcommand{\cE}{{\mathcal E}}

\newcommand{\Ker}{\mtype{\rm Ker\,}}
\newcommand{\im}{\mtype{Im\,}}
\newcommand{\ord}{\mbox{ord}}

\newcommand{\mtype}{\operatorname}

\title[The LFED and LNED Conjectures for $\cA_n$] 
{The LFED and LNED Conjectures for Laurent Polynomial Algebras}

  \author{Wenhua Zhao}      
    \date{\today}
\address{Department of Mathematics, Illinois State University, Normal, IL 61761. Email: wzhao@ilstu.edu}

\begin{document}

\begin{abstract}
Let $R$ be an integral domain of characteristic zero, $x=(x_1, x_2, \dots, x_n)$ 
$n$ commutative free variables, and $\cA_n\!:=R[x^{-1}, x]$, i.e., the Laurent 
polynomial algebra in $x$ over $R$. In this paper we first classify all locally finite or locally nilpotent $R$-derivations and $R$-$\cE$-derivations of $\cA_n$, 
where by an $R$-$\cE$-derivation of $\cA_n$ we mean an $R$-linear map of the form 
$\operatorname{Id}_{\cA_n}-\phi$ for some $R$-algebra endomorphism $\phi$ of $\cA_n$.
In particular, we show that $\cA_n$ has no nonzero locally nilpotent 
$R$-derivations or $R$-$\cE$-derivations. Consequently, the LNED conjecture proposed in \cite{Open-LFNED} for $\cA_n$ follows. 
We then show some cases of the LFED conjecture proposed in \cite{Open-LFNED} 
for $\cA_n$. In particular, we show that both the LFED and LNED conjectures hold for the Laurent polynomial algebras in one or two commutative free variables over a field of characteristic zero. 
\end{abstract}

   








%
%



\thanks{The author has been partially supported 
by the Simons Foundation grant 278638}

 \bibliographystyle{alpha}
    \maketitle


\renewcommand{\theequation}{\thesection.\arabic{equation}}
\renewcommand{\therema}{\thesection.\arabic{rema}}
\setcounter{equation}{0}
\setcounter{rema}{0}
\setcounter{section}{0}

\section{\bf Introduction} \label{S1}

Let $R$ be a unital commutative ring and $\cA$ an $R$-algebra. 
We denote by $1_\cA$ or simply $1$ the identity element of $\cA$, if $\cA$ is unital, and 
$\I_\cA$ or simply $\I$ the identity map of $\cA$, if $\cA$ is clear in the context.  

An $R$-linear endomorphism $\eta$ of $\cA$ is said to be {\it locally nilpotent} (LN) 
if for each $a\in \cA$ there exists $m\ge 1$ such that $\eta^m(a)=0$, 
and {\it locally finite} (LF) if for each $a\in \cA$ the $R$-submodule spanned 
by $\eta^i(a)$ $(i\ge 0)$ over $R$ is finitely generated.    
For each $R$-linear endomorphism $\eta$ of $\cA$ we denote by $\im \eta$ 
the {\it image} of $\eta$, i.e., $\im\eta\!:=\eta(\cA)$, 
and $\Ker \eta$ the {\it kernel} of $\eta$.

By an $R$-derivation $D$ of $\cA$ we mean an $R$-linear map 
$D:\cA \to \cA$ that satisfies $D(ab)=D(a)b+aD(b)$ for all $a, b\in \cA$. 
By an $R$-$\cE$-derivation $\delta$ of $\cA$ 
we mean an $R$-linear map $\delta:\cA \to \cA$ such that for all 
$a, b\in \cA$ the following equation holds:
\begin{align}\label{ProdRule2}
\delta(ab)=\delta(a)b+a\delta(b)-\delta(a)\delta(b).  
\end{align}

It is easy to verify that $\delta$ is an $R$-$\cE$-derivation of $\cA$, 
if and only if $\delta=\I-\phi$ for some $R$-algebra endomorphism 
$\phi$ of $\cA$. Therefore an $R$-$\cE$-derivation is  
a special so-called $(s_1, s_2)$-derivation introduced 
by N. Jacobson \cite{J} and also a special 
semi-derivation introduced 
by J. Bergen in \cite{Bergen}. $R$-$\cE$-derivations have also been 
studied by many others under some different names such as 
$f$-derivations in \cite{E0, E} and 
$\phi$-derivations in \cite{BFF, BV}, etc..

Now let us recall the following notion of associative algebras 
introduced in \cite{GIC, MS}. Since all algebras in this paper  
are commutative, here we recall only the case for commutative 
algebras.

\begin{defi} \label{Def-MS}
Let $\cA$ be a commutative $R$-algebra. An $R$-subspace  
$V$ of $\cA$ is said to be a {\it Mathieu subspace} (MS)  
of $\cA$ if for all $a, b\in \cA$ with $a^m\in V$ for all $m\ge 1$, 
we have $a^mb \in V$ for all $m\gg 0$. 
\end{defi}
  
Note that a MS is also called a {\it Mathieu-Zhao space} 
in the literature (e.g., see \cite{DEZ, EN, EH}, etc.), 
as first suggested by A. van den Essen \cite{E2}. 
 
The introduction of this new notion  
is mainly motivated by the study in \cite{Ma, IC} of 
the well-known Jacobian conjecture (see \cite{K, BCW, E}). 
See also \cite{DEZ}. But, a more interesting aspect 
of the notion is that it provides a natural   
generalization of the notion of ideals.

Next, we recall the cases of the so-called LFED and LNED conjectures 
proposed in \cite{Open-LFNED} for commutative algebras. 
For the study of some other cases of these two 
conjectures, see \cite{EWZ}, \cite{Open-LFNED}--\cite{OneVariableCase}.

\begin{conj}\label{LFED-Conj}
Let $K$ be a field of characteristic zero, $\cA$ a commutative $K$-algebra 
and $\delta$ a LF (locally finite) $K$-derivation or 
$K$-$\cE$-derivation of $\cA$. Then the image 
$\im\delta$ of $\delta$ is a MS of $\cA$.   
\end{conj}

\begin{conj}\label{LNED-Conj}
Let $K$ be a field of characteristic zero, $\cA$ a commutative 
$K$-algebra and $\delta$ a LN (locally nilpotent) $K$-derivation or 
$K$-$\cE$-derivation of $\cA$. Then $\delta$ maps every 
ideal of $\cA$ to a MS of $\cA$.  
\end{conj}

Throughout the paper we refer Conjecture \ref{LFED-Conj} as 
the (commutative) LFED conjecture, and Conjecture \ref{LNED-Conj} 
the (commutative) LNED conjecture.  
 
In this paper we study some cases of the two conjectures above 
for the Laurent polynomial algebra $\cA_n\!:=R[x^{-1}, x]$ in $n$ 
commutative free variables $x=(x_1, x_2, \dots, x_n)$ over 
an integral domain $R$ of characteristic zero.

In section \ref{S2} we classify all LF or LN 
$R$-derivations of $\cA_n$ (see Propositions \ref{ClassiD} and \ref{ClassiHD}).
In particular, we show that $\cA_n$ has no nonzero LN 
$R$-derivations or $R$-$\cE$-derivations. Hence the LNED 
conjecture \ref{LNED-Conj} holds (trivially) for $\cA_n$.

In section \ref{S3} we show first that the LFED 
conjecture \ref{LFED-Conj} holds for all LF $R$-derivations 
of $\cA_n$ (See Theorem \ref{LNED-Case}). 
We also show the LFED conjecture \ref{LFED-Conj} with Condition (C1) on page \pageref{C1} on $R$ for some LF $R$-$\cE$-derivations of $\cA_n$ (See Lemmas \ref{Lma5.2.1}, \ref{URootCase}  and \ref{W0-Lma} and Proposition \ref{A1=pmI}). In particular, 
we show that both the LFED and LNED conjectures hold for the Laurent polynomial algebras in one or two variables over a field of characteristic zero 
(see Theorem \ref{One-Two-VariableCase}). For the other LF $R$-$\cE$-derivations $\delta$ of $\cA_n$, we give an explicit description for $\im \delta$ (see Theorem \ref{ImDesp}) and reduce the LFED conjecture \ref{LFED-Conj} for $\delta$ to a special case of a conjecture (see Conjecture \ref{Comm-GeneralDK-Conj}) proposed in \cite{Open-LFNED}.


\renewcommand{\theequation}{\thesection.\arabic{equation}}
\renewcommand{\therema}{\thesection.\arabic{rema}}
\setcounter{equation}{0}
\setcounter{rema}{0}

\section{\bf Classification of Locally Finite or Locally Nilpotent 
Derivations and $\cE$-Derivations of $\cA_n$}\label{S2}

Throughout this paper $R$ stands for an integral domain of characteristic zero 
and $x=(x_1, x_2, \dots, x_n)$ $n$ commutative free variables. 
We denote by $R[x^{-1}, x]$ or $\cA_n$ the algebra of the Laurent polynomials 
in $x$ over $R$. The main purpose of this section is to classify all LF 
(locally finite) or LN (locally nilpotent) $R$-derivations and 
$R$-$\cE$-derivations of $\cA_n$. 

First, we fix the following notations that will be used in the rest of this paper. The notations introduced in the previous section will also be freely used.\\

{\bf Notations and Conventions:}
\vspace{3mm}
\begin{enumerate}
\item[$i)$] We denote by $R^\ast$ the set of units of $R$. 
For each $k\ge 1$ we let $\{e_i\,| \, 1\le i\le k\}$ be the 
standard basis of $\bZ^k$, $M_k(\bZ)$ the set of all $k\times k$ matrices 
with integer entries, and $\I_k$ the $k\times k$ identity matrix.

\item[$ii)$] For each nonzero $f\in \cA_n$ we denote by $\deg f$ the (generalized) {\it degree} of $f(x)$, and $\ord f$ the {\it order} of $f$, i.e., the minimum of all the degrees of the monomials in $f(x)$ 
with nonzero coefficients. Furthermore, we denote by $\supp f$ the {\it support} 
of $f$, i.e., the set of all $\alpha\in \bZ^n$ such that the coefficient of $x^\alpha$ is not zero, and $\poly f$ the {\it polytope} of $f$, i.e., the convex subset of $\bR^n$ spanned by $\supp f$.   

\item[$iii)$] We denote by $\p_i$ $(1\le i\le n)$ the $R$-derivation $\p/\p x_i$ of $\cA_n$. Then every $R$-derivation $D$ of $\cA_n$ can be written uniquely as $\sum_{i=1}^n a_i(x)$ for some $a_i(x)\in \cA_n$. For all $k\ge 0$, we let 
$\cH_k$ denote the $R$-subspace of $\cA_n$
of homogeneous Laurent polynomials 
of the (generalized) degree $k$. \\
\end{enumerate}

We start with a classification of all LF  or LN $R$-derivations of $\cA_n$. 

\begin{propo}\label{ClassiD}
Let $D$ be an $R$-derivation of $\cA_n$.   Then
\begin{enumerate}
\item[$1)$] $D$ is locally finite, if and only if  
$D=\sum_{1\le i\le n} q_ix_i\p_i$ for some $q_i\in R$.
  \item[$2)$] $D$ is locally nilpotent, (if and) only if $D=0$, i.e., $\cA$ has no
nonzero locally nilpotent $R$-derivations. 
\end{enumerate}
\end{propo}

\pf $1)$ The $(\Leftarrow)$ part is obvious. To show the $(\Rightarrow)$ part
we may obviously assume $D\ne 0$ and write $D$ uniquely as $D=\sum_{k=r}^s D_k$ 
with $r\le s$ and $D_k$ $(r\le k\le s)$ an $R$-derivation of $\cA_n$ such that 
$D_k\cH_i\subseteq\cH_{i+k}$ for all $i\in \bZ$ and $D_r, D_s \ne 0$. 

We first show $r=s=0$, i.e., $D=D_0$. Below we give only a proof 
for $s=0$ by using the degree of elements of $\cA_n$. The proof 
for $r=0$ is similar (by using the order of elements of $\cA_n$, 
instead). 

Assume $s\ne 0$. Then for all $1\le i\le n$ and $m\ge 1$, 
$D_s^m x_i$ and $D_s^m x_i^{-1}$, if not zero, 
are homogeneous of the (generalized) degree $ms+1$ and $ms-1$, 
respectively. 
On the other hand, since $D$ is LF,  
$D^m x_i$ and $D^m x_i^{-1}$ $(1\le i\le n$ and $m\ge 1)$ lie 
inside a finitely generated $R$-submodule of $\cA_n$. 
Note that $D_s^m x_i$ and $D_s^m x_i^{-1}$, if not zero, 
are the leading terms of $D^m x_i$ and $D^m x_i^{-1}$, 
respectively. Therefore, there exists $N\ge 1$ such that  
$D_s^m x_i=D_s^m x_i^{-1}=0$ for all 
$1\le i\le n$ and $m\ge N$. 

Since $D_s x_i^{-1}=-x_i^2D_s x_i$ for each $1\le i\le n$, 
we have that $D_sx_i=0$, if and only if $D_s x_i^{-1}=0$. Since $D_s\ne 0$, 
there exists $1\le i\le n$ such that $D_s x_i, D_s x_i^{-1}\ne 0$.
Let $k$ (resp., $\ell$) be the least positive 
integer such that $D_s^k x_i=0$ (resp., $D_s^\ell x_i^{-1}=0)$. 
  
If $k+\ell\ge 3$, then by the Leibniz rule we have 
$$
0=D_s^{k+\ell-2} 1=  D_s^{k+\ell-2}(x_ix_i^{-1})=\binom{k+\ell-2}{k-1} 
(D_s^{k-1}x_i)(D_s^{\ell-1}x_i^{-1}).
$$ 
Since $\cA_n$ is an integral domain of characteristic zero, we have 
$D_s^{k-1}x_i=0$ or $D_s^{\ell-1}x_i^{-1}=0$. But this contradicts to 
either the choice of $k$ or that of $\ell$. Therefore 
$k+\ell \le 2$, whence $k=\ell=1$, i.e., $D_s x_i=D_sx_i^{-1}=0$.
Contradiction. Hence $s=0$, as desired. 
 
Now we may write $D=D_0=\sum_{i=1}^n a_i(x)\p_i$ for some homogeneous 
$a_i(x)\in \cA_n$ $(1\le i\le n)$ of degree one. 
We need to show that for all $1\le i\le n$ 
we have $a_i(x)=q_ix_i$ for some $q_i\in R$.
But, up to the conjugation by a permutation (on $x_i$'s) automorphism of $\cA_n$ 
it suffices to show only the case for $a_1(x)$. 

Write $a_1(x)=b(x)\p_1+q_1x_1\p_1$ with $q_1\in R$ and $b(x)$ 
homogeneous of degree one and independent on $x_1$. 
Then for all $m\ge 1$, it is easy to see that 
$D^m x_1^{-1}$ as a Laurent polynomial in $x_1$ 
over $R[x_i^{-1}, x_i\,|\, 2\le i\le n]$ 
has the least degree (in $x_1$) 
term $(-1)^m m! b^m(x)x_1^{-m-1}$. 
Since $D$ is LF (again), we have $b(x)^mx_1^{-m-1}=0$ 
for all $m\gg 0$, whence $b(x)=0$ and 
statement $2)$ follows. 
 
$2)$ Let $D$ be a LN (locally nilpotent) $R$-derivation of $\cA_n$. 
Then $D$ is LF and by statement $1)$ we have 
$D=D_0=\sum_{1\le i\le n} q_ix_i\p_i$ 
for some $q_i\in R$. If $D\ne 0$, then $q_j\ne 0$ for some $1\le j\le n$, 
and for all $m\ge 1$ we have $D^m x_j^{-1}=(-1)^m q_j^m x_j^{-1}\ne 0$. 
Contradiction.
\epfv
 
\begin{rmk}
By some similar arguments as that in the proof of Proposition \ref{ClassiD} it is easy to see that Proposition \ref{ClassiD} also holds for the Laurent polynomial algebra in noncommutative free variables over $R$. 
\end{rmk}
 
Next, we give a classifications of all LF or LN   
$R$-$\cE$-derivations of $\cA_n$. 
Let $\phi\ne 0$ be an $R$-algebra endomorphism 
of $\cA_n$.  
Since $\phi(1)$ is an idempotent of $\cA_n$ and $\cA_n$ is an integral domain, 
we have $\phi(1)=0$ or $1$. Since $\phi(1)=0$ implies $\phi=0$, we have $\phi(1)=1$, 
i.e., $\phi$ preserves the unity $1$ of $\cA_n$. Consequently, 
it also preserves the units of $\cA_n$. 

Since all units of $\cA_n$ have the form $qx^\alpha$ with $q\in R^*$ 
and $\alpha\in \bZ^n$, 
we see that $\phi(x_i)=q_i x^{\alpha_i}$ with $q_i\in R^*$ and 
$\alpha_i\in \bZ^n$ for all $1\le i\le n$. Set $q\!:=(q_1, q_2, \dots, q_n)$ 
and for each $1\le i\le n$ write $\alpha_i=(\alpha_{1i}, \alpha_{2i}, \dots, \alpha_{ni})$ and form the $n\times n$ matrix $A\!:=(\alpha_{ij})\in M_n(\bZ)$. 
For convenience, for the case $\phi=0$ we simply let $A=0$ and 
$q_i=0$ $(1\le i\le n)$. 
For the rest of this paper we call the matrix $A$ the 
{\it exponent matrix} of the $R$-algebra endomorphism 
of $\phi$.

With the remarks and notations above it is easy to see that    
$\phi(x^\beta)=q^{\beta}x^{A\beta}$ for all $\beta\in \bZ^n$. 
In particular, both $\phi$ and $\I-\phi$ maps each monomial to a scalar multiple of a monomial. Therefore, $\I-\phi$ is LF, if and only if $\phi$ is LF, if and only if  for each $\alpha\in \bZ$ the set 
$\{\phi^m(x^\alpha)\,|\, m\ge 0\}$ is finite.

\begin{propo}\label{ClassiHD}
With the notations above, for all $R$-endomorphism $\phi$ of $\cA_n$ 
the following statements hold:  
\begin{enumerate}
\item[$1)$] $\I-\phi$ is locally finite if and only if there exist $0\le k\le n$, an invertible and finite order $B \in M_{n-k}(\bZ)$, and an invertible 
$S\in M_n(\bZ)$ such that 
\begin{align}\label{ClassiHD-eq1}
A=S
\begin{pmatrix}
0&0\\
0&B
\end{pmatrix}
S^{-1}.
\end{align} 
 \item[$2)$] $\I-\phi$ is locally nilpotent (if and) only if $\phi=\I$, i.e., $\cA_n$ has no nonzero locally nilpotent $R$-$\cE$-derivations. 
\end{enumerate}
\end{propo}
 
\pf $1)$ The $(\Leftarrow)$ part can be checked easily. 
To show the $(\Rightarrow)$ part, we assume $A\ne 0$, 
for the case $A=0$ is trivial, 
and set for all $m\ge 1$ 
\begin{align}\label{ClassiHD-peq0}
\tilde A_m\!:=\I_{n\times n}+A+A^2+\cdots+A^{m-1}.
\end{align} 
Then it is easy to see inductively that for all $\alpha\in \bZ^n$ and $m\ge 1$ we have
\begin{align}\label{ClassiHD-peq2}
\phi^m(x^\alpha)=q^{\tilde A_m\alpha}x^{A^m\alpha}.
\end{align}

Since $\I-\phi$ is LF, $\phi$ is also LF. 
Then for each $1\le i\le n$ the set $\{\phi^m(x_i)\,|\, m\ge 1\}$ 
is finite, whence by Eq.\,(\ref{ClassiHD-peq2}) 
there exists $1\le k_i<\ell_i$ such that 
\begin{align}\label{ClassiHD-peq2.1}
A^{k_i}e_i=A^{\ell_i}e_i,
\end{align} 
where $e=(e_1, e_2, \dots, e_n)$ is the standard basis of $\bZ^n$. 

Multiplying a power of $A$ to Eq.\,(\ref{ClassiHD-peq2.1}) 
for each $1\le i\le n$ we may assume $k_i=k_j$ for 
all $1\le i, j\le n$. We denote this integer by $k$ and 
set $\ell=k+\prod_{i=1}^n(\ell_i-k)$. Then 
it is easy to see that $A^k e_i=A^\ell e_i$ 
for all $1\le i\le n$, whence 
$A^k=A^\ell$. 


Let $H$ be the subgroup of $\bZ^n$ formed by $\alpha\in \bZ^n$ 
such that $A\alpha=0$. If $H=0$, then $A$ is of rank $n$, 
and from the fact $A^k=A^\ell$ with $1\le k<\ell$, 
we see that $A^{k-\ell}=\I_n$. Hence the statement 
holds in this case. So we assume $H\ne 0$.

By a well-known fact (e.g., see \cite[Theorem 1.6]{Hun}) on subgroups 
of free abelian groups  
there exists a basis $(\beta_1, \beta_2, \dots, \beta_n)$ of 
$\bZ^n$ such that $H$ is the (free) subgroup generated by 
$(d_1\beta_1, d_2\beta_2, \dots, d_k\beta_k)$ for some  
$1\le k\le n$ and positive integers $d_i$ $(1\le i\le k)$. 
Since $A(d_i\beta_i)=0$ implies $A\beta_i=0$ for each  
$1\le i\le k$, we hence have $d_i=1$ for 
all $1\le i\le k$, and also $k\le n-1$, for $A\ne 0$. 
Therefore there exist $k\times (n-k)$ matrix $C$, and  
$B\in M_{n-k}(\bZ)$ of rank $n-k$, and an invertible $S\in M_n(\bZ)$ 
such that  
\begin{align}\label{ClassiHD-peq1}
A=S
\begin{pmatrix}
0&C\\
0&B
\end{pmatrix}
S^{-1}.
\end{align}

Set 
$D\!:=
\begin{pmatrix} 
0&C\\
0&B
\end{pmatrix}
$.
Since   
$D^m=
\begin{pmatrix} 
0&CB^{m-1}\\
0&B^m.
\end{pmatrix}
$
for all $m\ge 1$, by the fact $A^k=A^\ell$ 
proved above we have $B^k=B^\ell$, 
which implies $B^{k-\ell}=\I_{n-k}$. Therefore $B$ is invertible and 
of finite order. So we need only to show $C=0$. 

Assume otherwise, i.e., $C\ne 0$. 
Denote by $\psi$ the $R$-algebra automorphism of $\cA_n$ that maps 
$x^\alpha$ to $x^{S\alpha}$ for all $\alpha\in \bZ^n$. 
Replacing $\phi$ by $\psi^{-1}\phi\psi$ we may assume 
$A=D$. Since $\phi$ is LF, so is $\phi^m$ for all $m\ge 1$.
Replacing $\phi$ by $\phi^{k-\ell}$ we may assume $B=\I_{n-k}$, i.e., 
$A= \begin{pmatrix} 
0&C'\\
0&\I_{n-k}
\end{pmatrix}
$,
where $C'=CB^{j-i-1}\ne 0$.

Now let $y=(x_1, x_2, \dots, x_k)$, $z=(x_{k+1}, x_{k+1}, \dots, x_n)$, 
and $\beta\in \bZ^{n-k}$ such that $C\beta\ne 0$. Then for all $m\ge 1$,  
$\phi^m(z^\beta)$ up to an invertible multiplicative scalar is equal to 
$y^{mC'\beta}z^\beta$. Since $mC'\beta$ $(m\ge 1)$ are all distinct, 
$\im \phi$ contains infinitely distinct 
monomials, which contradicts to the assumption 
that $\phi$ is LF.  Therefore $C'=0$, and hence $C=0$, as desired. 

$2)$ Assume that $\I-\phi$ is LN. Then $\I-\phi$ 
is also LF, and statement $1)$ holds for $\I-\phi$. 
In particular, there are exist some $1\le k<\ell$ such that $A^k=A^\ell$.  

We first assume $A=A^2$. Then for each $\alpha\in \bZ^n$ 
there exists $N\ge 1$ such that $0=(\I-\phi)^N x^\alpha=x^\alpha+cx^{A\alpha}$ 
for some $c\in R$, whence $A\alpha=\alpha$. Therefore $A=\I_n$. 
Consequently, for all $\alpha\in \bZ^n$ 
and $m\ge 1$, we have 
$\phi(x^\alpha)=q^\alpha x^{\alpha}$ and $(\I-\phi)^m x^\alpha=(1-q^\alpha)^m x^\alpha$. 
Since $\I-\phi$ is LN, we have $q^\alpha=1$ for all $\alpha\in \bZ$, 
whence $q_i=1$ for all $1\le i\le n$ and $\phi=\I$, as desired.  

For the general case, let $r\ge 1$ be the order of the matrix $B$ in 
Eq.\,(\ref{ClassiHD-eq1}). Then by Eq.\,(\ref{ClassiHD-eq1}) we have $A^r=(A^r)^2$,  
and  by Eq.\,(\ref{ClassiHD-peq2}) 
$\phi^r(x^\alpha)=  q^{\tilde A_r\alpha} x^{A^r\alpha}$ for all $\alpha\in \bZ^n$.
Since $(\I-\phi^r)=(\sum_{i=1}^{r-1}\phi^i)(\I-\phi)$, 
$(\I-\phi^r)$ is also LN. Applying the case $A=A^2$ shown above to 
the $R$-algebra endomorphism $\phi^r$ we get $\phi^r=\I$. 
Then by \cite[Corollary 6.5]{Open-LFNED} we have $\phi=I$.

A more straightforward way to show the last step above is as follows. 

For each $\alpha\in \bZ$, we have $(\I-\phi^r)x^\alpha=0=(\I-\phi)^N x^\alpha$ 
for some $N\ge 1$. On the other hand, both the polynomials $(1-t^r)$ and $(1-t)^N$ 
are monic and have the ``greatest common divisor" $1-t$ in $\bC[t]$ . 
By the Euclidean algorithm it is easy to verify that there exist 
$u(t), v(t)\in \bZ[t]\subseteq R[t]$ such that 
$(1-t^r)u(t)+(1-t)^Nv(t)=1-t$. Substituting $\phi$ for $t$ in the equation above 
and then applying it to $x^\alpha$ we see that $(\I-\phi)(x^\alpha)=0$
for each $\alpha\in \bZ$.  Hence $\phi=\I$, as desired. 
\epfv


\renewcommand{\theequation}{\thesection.\arabic{equation}}
\renewcommand{\therema}{\thesection.\arabic{rema}}
\setcounter{equation}{0}
\setcounter{rema}{0}

\section{\bf The LFED and LNED Conjectures for $\cA_n$}\label{S3}

In this section we consider the LFED conjecture \ref{LFED-Conj} 
and the LNED conjecture \ref{LNED-Conj} for the Laurent 
polynomial algebras over an integral domain 
of characteristic zero. 

Throughout this section $R$ denotes an integral domain 
$R$ of characteristic zero, $x=(x_1, x_2, \dots, x_n)$ commutative 
free variables and $\cA_n=R[x^{-1}, x]$. 
All other notations and conventions introduced 
in the previous two sections will also be in force 
in this section. 

First, we have the following:  

\begin{theo}\label{LNED-Case}
$1)$ The LNED conjecture \ref{LNED-Conj} holds for $\cA_n$.

$2)$ The LFED conjecture \ref{LFED-Conj} holds for all 
LF $R$-derivations of $\cA_n$.
\end{theo}

\pf $1)$ follows trivially from Propositions \ref{ClassiD} and \ref{ClassiHD}, 
by which $\cA_n$ has no nonzero LN $R$-derivations or LN $R$-$\cE$-derivations.
$2)$ follows directly from Propositions \ref{ClassiD} and the lemma below.
\epfv

\begin{lemma}
Let $D=\sum_{1\le i\le n} a_ix_i\p_i$ for some $a_i\in R$.   
Then $\im D\!:=D(\cA_n)$ is a MS of $\cA_n$. 
\end{lemma}

The lemma for polynomial algebras $K[x]$ over a field $K$ of characteristic zero 
has been proved in Lemma $3.4$ in \cite{EWZ}. 
It is easy to see that the proof there with some 
slight modifications also works for this more  
general case. So we skip the proof here. 
 
Next, we consider the $\cE$-derivation case of the LFED 
conjecture \ref{LFED-Conj} for $\cA_n$. We fix a LF $R$-algebra endomorphism 
$\phi$ of $\cA_n$. Then by Proposition \ref{ClassiHD} and up to a 
conjugation of $\phi$ we may assume that $\phi$ maps 
$x^\alpha$ $(\alpha\in \bZ^n)$ to $q^\alpha x^{A\alpha}$, where  
$q=(q_1, q_2, \dots, q_n)$ with $q_i\in R^*$ $(1\le i\le n)$ 
and $A\in M_n(\bZ)$ is the {\it exponent matrix} of $\phi$ 
given by  
\begin{align} \label{Form4A}
A= 
\begin{pmatrix}
0&0\\
0&B
\end{pmatrix}
\end{align} 
for some invertible $B\in M_{n-k}(\bZ)$ of finite order.

We start with the following simple case.

\begin{lemma}\label{Lma5.2.1}
Assume $A=0$. Then $\im(\I-\phi)=\Ker\phi$, which is the maximal ideal 
of $\cA_n$ generated by $(x_i-q_i)$ $(1\le i\le n)$.  
\end{lemma}
 
\pf Note that in this case $\phi(x^\alpha)=q^\alpha$ for all $\alpha\in \bZ^n$. In particular, 
$\phi^2=\phi$. Then it is easy to check directly 
(or by the more general 
\cite[Proposition 5.2]{Open-LFNED}) 
that the lemma indeed holds.
\epfv
  
Next we drive the following reduction. 

\begin{lemma}\label{FuRankRed}
To show whether or not $\im(\I-\phi)$ is a MS, we may assume that the exponent matrix 
$A$ of $\phi$ is invertible and of finite order. 
\end{lemma}

\pf Let $A$ be given as in Eq.\,(\ref{Form4A}). 
If $k=0$, then there is nothing to show. 
If $k=n$, i.e., $A=0$, then by Lemma \ref{Lma5.2.1} $\im(\I-\phi)$ is 
an ideal of $\cA$, hence a MS of $\cA$. So we assume $k\ne 1, n$.
 
%
%
%

Let $J$ be the ideal of $\cA_n$ generated by 
$x_i-q_i$ $(1\le i\le k)$. Then $\phi(J)=0$ and 
$(\I-\phi)(J)=J$, whence $J\subseteq \im(\I-\phi)$.
 
Let $\bar\cA_n=\cA_n/J$ and $\bar\phi$ be the $R$-algebra endomorphism of $\bar\cA_n$ 
induced by $\phi$. Then it can be readily verified that 
$\im(\I-\phi)/J=\im(\I_{\bar\cA_n}-\bar\phi)$. Then by \cite[Proposition 2.7]{MS}
we have that $\im(\I-\phi)$ is a MS of $\cA_n$, if and only if 
$\im(\I_{\bar\cA_n}-\bar\phi)$ is a MS of $\bar\cA_n$. 
Note that $\bar\cA$ as an $R$-algebra is isomorphic 
to the Laurent polynomial algebra in $x_j$ $(k+1\le j\le n)$ 
and the exponent matrix of $\bar\phi$ is $B\in M_{n-k}(\bZ)$, 
which is invertible and of finite order. Hence the lemma follows.
\epfv

From now on by Lemma \ref{FuRankRed} above we may and will assume that the exponent 
matrix $A$ of $\phi$ is invertible and of finite order. We denote by $r$ 
the {\it order} of $A$ and assume further 
\footnote{Without the condition (\ref{C1}) the LFED conjecture \ref{LFED-Conj} 
may not be true for $\I-\phi$, e.g., see \cite[Example 2.8]{Open-LFNED}.} 
\begin{align} \label{C1} 
\bQ(q)\subseteq R, \tag{C1} 
\end{align}
where $\bQ(q)$ is the field of fractions 
of $\bZ[q]\subseteq R$.

We also need to fix the following notations:
\begin{align}
\tilde A_m\!:&=\I_n+A+\cdots+A^{m-1} \,\,\text{ (for all } m\ge 1). \label{Def-Am}\\
\tilde A\!:&=\I_n+A+\cdots+A^{r-1}.\label{Def-Ar}\\
W \!:&=\{\alpha\in \bZ^n \,|\, q^{\tilde A \alpha}=1\}.\label{DefW}
\end{align} 

 
\begin{lemma}\label{URootCase}
If the abelian subgroup $W$ of $\bZ^n$ has rank $n$ 
(e.g., when $q_i$ $(1\le i\le n)$ are all roots of unity).  
Then $\im(\I-\phi)$ is a MS of $\cA_n$.
\end{lemma}

\pf Let $H=\bZ^n/W$. Then $H$ is a finite abelian group. Let 
$d$ be the order of $H$, i.e., $d=|H|$. Then $d\alpha\in W$ for all 
$\alpha\in \bZ^n$. Let $s=rd$. Then $A^s=\I_n$ and $\tilde A_s=d\tilde A$. 
By Eqs.\,(\ref{ClassiHD-peq2}) and (\ref{Def-Am})-(\ref{DefW}) 
we have for all $\alpha\in \bZ^n$  
$$
\phi^s(x^\alpha)=q^{\tilde A_s\alpha}x^{A^s\alpha}=q^{d\tilde A\alpha}x^\alpha.
$$
Therefore $\phi^s=\I$, and by \cite[Corollary 5.5]{Open-LFNED} the lemma follows.  
\epfv

Before we show the next lemma, let us recall the following 
remarkable Duistermaat-van der Kallen 
Theorem \cite{DK}.
 
\begin{theo}\label{DvK-Thm}
Let $K$ be a field of characteristic zero, $z=(z_1, \dots, z_n)$ 
commutative free variables and $M$ the $K$-subspace of $K[z^{-1}, z]$ 
of the Laurent polynomials with no constant term. Then for each nonzero 
$f\in K[z^{-1}, z]$ such that $f^m\in M$ for all $m\ge 1$ we have   
$0\not\in \poly f$.   
Consequently, $M$ is a MS of $K[z^{-1}, z]$.  
\end{theo}

\begin{lemma}\label{W0-Lma}
Assume $W=0$ (e.g., when $q^\alpha\ne 1$ for all $0\ne \alpha\in \bZ^n$). Then 
\begin{enumerate}
  \item[$1)$] $\im(\I-\phi)$ is the $R$-subspace spanned by 
all Laurent polynomial without constant terms.
\item[$2)$]  $\im(\I-\phi)$ is a MS of $\cA_n$.
\end{enumerate}
\end{lemma}

\pf $1)$ For all $\alpha\in \bZ^n$, by Eq.\,(\ref{ClassiHD-peq2}) we have 
$\phi^r(x^\alpha)=q^{\tilde A \alpha}x^\alpha$, whence 
$(\I-\phi^r)(x^{\alpha})=(1-q^{\tilde A \alpha})x^\alpha$.
Since $W=0$, $q^{\tilde A \alpha}\ne 1$ for all $0\ne \alpha\in \bZ^n$.
Then by Condition (\ref{C1}) $x^\alpha\in \im(\I-\phi^r)$ for all 
$0\ne \alpha\in \bZ^n$.

Furthermore, since $(\I-\phi^r)=(\I-\phi)\sum_{i=1}^{r-1}\phi^i$, 
we have $\im (\I-\phi^r)\subseteq \im(\I-\phi)$. Hence 
$x^\alpha\in \im(\I-\phi)$ for all $0\ne \alpha\in \bZ^n$.
Since $\im(\I-\phi)$ is a homogeneous $R$-subspace of $\cA_n$ 
and $1\not \in \im(\I-\phi)$, the statement follows.

$2)$ Let $K_R$ be the field of fractions of $R$ and 
$f\in \cA_n$ such that $f^m\in \im(\I-\phi)$ for all $m\ge 1$. 
Viewing $f$ as an element of $K_R[x^{-1}, x]$ and by Theorem \ref{DvK-Thm} 
we get $0\not\in \poly f$. Hence for all 
$g\in\cA_n$ we have $0\not \in \poly (f^mg)$ 
when $m\gg 0$.  Therefore $\im(\I-\phi)$ 
is a MS of $\cA_n$.
\epfv 

From now on we assume $W\ne 0$. By a well-known fact 
(e.g., see \cite[Theorem 1.6]{Hun}) on subgroups of free abelian groups  
there exists a basis  $(\alpha_1, \alpha_2, \cdots, \alpha_n)$ of $\bZ^n$ 
such that $W$ is freely spanned by $(d_1\alpha_1, d_2\alpha_2, \cdots, d_k\alpha_n)$ 
for some $1\le k\le n$ and positive integers $d_i$ with 
$d_1\,|\, d_{2}\,| \cdots\,|\, d_{k-1}\,|\, d_k$. 

By conjugating an automorphism of $\cA_n$ to $\phi$ 
we may further assume that the basis $(\alpha_1, \alpha_2, \cdots, \alpha_n)$ 
is the standard basis $(e_1, e_2, \cdots, e_n)$ of $\bZ^n$. 
Furthermore, we also fix the following notations for the rest of this section:
\begin{enumerate}  
\item[$i)$] Denote by $E_1$ (resp., $E_2$) the subgroup of $\bZ^n$ 
generated by $e_i$ with $1\le i\le k$ (resp., $k+1\le i\le n$). 
\item[$ii)$] Set $\ell=n-k$; $y=(y_1, y_2, \dots, y_k)$ with $y_i=x_i$ $(1\le i\le k)$; and 
$z=(z_1, z_2, \dots, z_\ell)$ with $z_j=x_{k+j}$ ($1\le j\le \ell$).  
\item[$iii)$] Set $\cB_1\!:=R[y_1, y_2,\dots, y_k]$ and $\cB_2\!:=R[z_1, z_2,\dots, z_\ell]$.
\item[$iv)$] Set $\phi_1\!:=\phi\,|\,_{{}_{\cB_1}}$ and $\phi_2\!:=\phi\,|\,_{{}_{\cB_2}}$.    
\end{enumerate}

\begin{lemma}\label{Lma(-2)}
With the notations fixed above the following statements hold:
\begin{enumerate} 
  \item[$1)$] for all $\alpha\in \bZ^n$, 
   $q^{\tilde A \alpha}$ is a root of unity, if and only if  $\alpha\in E_1$. 
  \item[$2)$] there exist $A_1\in M_k(\bZ)$, $A_2\in M_{n-k}(\bZ)$ and 
a $k\times {(n-k)}$ matrix $B$ with integer entries such that   
\begin{align}\label{Lma(-2)-eq1}
A_1^r=\I_k, \qquad A_2^r=\I_{n-k}, \qquad                      
A=\begin{pmatrix} 
A_1& B\\
0&A_2
\end{pmatrix}                     
\end{align}
\item[$3)$]  \quad  
  $\phi (\cB_1)\subseteq \cB_1$, and\, 
$\cB_1\cap \im(\I-\phi)=\im(\I_{\cB_1}-\phi_1)$  
\end{enumerate}
\end{lemma}

\pf $1)$ $(\Rightarrow)$ Since $ \alpha\in E_1$, $d_k\alpha\in W$. Then  
$1=q^{\tilde A(d_k\alpha)}=q^{d_k\tilde A \alpha}=(q^{\tilde A \alpha})^{d_k}$, 
whence $q^{\tilde A \alpha}$ is a root of unity.

To show the $(\Leftarrow)$ part, since $\bZ^n=E_1\oplus E_2$, it suffices 
to show that for any $0\ne \beta\in E_2$, $q^{\tilde A \beta}$ is not 
a root of unity.

Assume otherwise, say, $1=(q^{\tilde A \beta})^j=q^{\tilde A (j\beta)}$ for some $j\ge 1$. Then $j\beta\in W\subseteq E_1$. Hence $j\beta \in E_1 \cap E_2=0$ and $\beta=0$.  Contradiction.

$2)$ Since $A^r=\I_n$, 
we have $\tilde A A=A\tilde A=\tilde A$. 
For each $1\le i\le k$, since $d_ie_i\in W$, we have 
$$
q^{\tilde A (A d_i e_i)}=q^{(\tilde A A) (d_i e_i)}= q^{\tilde A (d_i e_i)}=1.
$$
Therefore $d_i (A e_i)\in W$, whence $Ae_i\in E_1$, 
whence $AE_1\subseteq E_1$. Consequently, $A$ has the form as the matrix 
in Eq.\,(\ref{Lma(-2)-eq1}). The first two equations in Eq.\,(\ref{Lma(-2)-eq1}) 
follow immediately from the last one and the fact $A^r=\I_n$.

$3)$ By statement $2)$ we immediately have $\phi(\cB_1)\subseteq\cB_1$ and $(\I-\phi)(\cB_1)\subseteq\cB_1$.  Hence $\im(\I_{\cB_1}-\phi_1) \subseteq \cB_1\cap \im(\I-\phi)$.
Conversely, let $b(y)\in \cB_1\cap \im(\I-\phi)$. 
Then there exists $f(x)\in \cA_n$ such that $(\I-\phi)(f)=b(y)$. 

Write $f(x)=a_0(y)+\sum_{0\ne \beta\in E_2} a_\beta(y) z^{\beta_i}$ 
with $a_\beta(y)\in \cB_1$ $(\beta\in E_2)$. 
Note that by statement $2)$ the projection of  
$A\beta$ on $E_2$ for each nonzero $\beta\in E_2$ is also nonzero, 
since $A_2$ in Eq.\,(\ref{Lma(-2)-eq1}) is invertible.   
Hence $(\I-\phi)$ preserves the $R$-subspace 
$\sum_{0\ne \beta\in E_2} z^\beta\cB_1$ of $\cA_n$, by which  
the term of $(\I-\phi)(f)$ that lies in $\cB_1$ is equal 
to $(\I-\phi)(a_0)$, whence $(\I-\phi)(a_0)=b(y)$.  
Therefore $b(y)\in \im(\I_{\cB_1}-\phi_1)$, and 
the statement follows. 
\epfv

%
%
 


\begin{lemma}\label{Lma(-1)}
Set $s\!:=rd_k$, $\tilde q_i\!:=q^{\tilde A_s e_{i+k}}$ for all $1\le i\le \ell=n-k$, 
and $\tilde q\!:=(\tilde q_1, \tilde q_2, \dots, \tilde q_\ell)$. 
Then the following statements hold:
\begin{enumerate}
  \item[$1)$] $\phi^s(z^\beta)=\tilde q^\beta z^\beta$ for all $\beta\in E_2$.
  In particular, $\phi^s(\cB_2)\subseteq \cB_2$.   
 \item[$2)$] For all $0\ne \beta\in E_2$, $\tilde q^{\beta}$ is not a root of unity. 
 \item[$3)$] $\phi^s\,|_{\cB_1}=\I_{\cB_1}$. 
\end{enumerate}
\end{lemma}

\pf $1)$ 
 Since $r\,|\, s$, we have $A^s=\I_n$ and 
$\tilde A_s=d_k\tilde A$. Then by Eq.\,(\ref{ClassiHD-peq2})  we have for all 
$\alpha\in \bZ^n$    
\begin{align}\label{Lma(-1)-peq1}
\phi^s(x^\alpha)=q^{\tilde A_s\alpha}x^\alpha=
q^{d_k\tilde A\alpha}x^\alpha, 
\end{align}
from which by letting 
$\alpha=\beta\in E_2$ statement $1)$ follows  .

$2)$ Assume otherwise, then $q^{\tilde A\beta}$ 
is also a root of unity, since  
$\tilde q^{\beta}=q^{\tilde A_s\beta}=q^{d_k\tilde A\beta}=
(q^{\tilde A\beta})^{d_k}$. By Lemma \ref{Lma(-2)}, $1)$ 
we have $\beta\in E_1\cap E_2$, whence $\beta=0$. Contradiction.

$3)$ By Lemma \ref{Lma(-2)}, $2)$ we have $AE_1\subseteq E_1$, 
whence $\tilde A_s E_1 \subseteq E_1$. Since $d_k E_1\subseteq W$, 
by Eq.\,(\ref{Lma(-1)-peq1}) the statement follows.  
\epfv

Now we give a complete description for $\im(\I-\phi)$.  

\begin{theo}\label{ImDesp}
$\im(\I-\phi)$ is the $R$-subspace of $\cA_n$ consisting of all 
$f(x)\in\cA_n$ of the form   
\begin{align}\label{ImDesp-eq1}
f(x)=a_0(y)+\sum_{0\ne \beta\in E_2} 
 a_\beta(y) z^\beta
\end{align}
with $a_\beta(y)\in \cB_1$ $(0\ne \beta\in E_2)$ and 
$a_0(y)\in \im(\I_{\cB_1}-\phi_1)$. 
\end{theo}

\pf Let $s=rd_k$ and $\tilde q$ as in Lemma \ref{Lma(-1)}, 
by which we have $\phi^s\,|_{\cB_1}=\I_{\cB_1}$. So we may view $\phi^s$ 
as a $\cB_1$-algebra endomorphism of $\cA_n=\cB_1[z^{-1}, z]$ such that  
$\phi^s(z^\beta)=\tilde q^\beta z^\beta$ for all $\beta\in E_2$.
Since $\tilde q^\beta$ is not a root of unity for all 
$0\ne \beta \in E_2$ and $\bQ[\tilde q]\subseteq\bQ[q]\subseteq R\subseteq \cB_1$, 
by Lemma \ref{W0-Lma} we have that $z^\beta\cB_1\subseteq\im(\I-\phi^s)$ 
for all $0\ne \beta\in E_2$. Since $(1-t)\,|\, (1-t^s)$, 
we have $\im(\I-\phi^s)\subseteq\im(\I-\phi)$, whence 
$z^\beta\cB_1\subseteq\im(\I-\phi)$ 
for all $0\ne \beta\in E_2$.

Now let $f(x)\in \cA_n$ and write $f(x)=a_0(y)+\sum_{0\ne \beta\in E_2} 
a_\beta(y) z^\beta$ with $a_\beta(y)\in \cB_1$ $(\beta\in E_2$. 
Then from the fact shown above we have that $f\in \im(\I-\phi)$, 
if and only if $a_0(y)\in\im(\I-\phi)$.   
Combining with Lemma \ref{Lma(-2)}, $3)$ we further have that 
$f\in \im(\I-\phi)$, iff $a_0(y)\in \im (\I-\phi)(\cB_1)
=\im(\I_{\cB_1}-\phi_1)$. Hence the theorem 
follows.
\epfv

%
%
 
One remark on Theorem above and the LFED conjecture  \ref{LFED-Conj} 
for $\cA_n$ is as follows.

Since $q_i$ $(1\le i\le k)$ are all roots of unity, by Lemma \ref{URootCase} we see that $\im(\I_{\cB_1}-\phi_1)$ is a MS of $\cB_1$. Another way to see this fact 
is to apply Lemma \ref{Lma(-1)}, $3)$  and \cite[Corollary 5.5]{Open-LFNED}. 
Then by Proposition \ref{ClassiHD}, $1)$ and Theorem \ref{ImDesp} we see that 
when $R$ is a field $K$ of characteristic zero, 
the LFED conjecture  \ref{LFED-Conj} for $\cA_n$ follows from the following conjecture, which is the case of \cite[Conjecture 4.4]{Open-LFNED} 
for commutative $K$-algebras.

\begin{conj}\label{Comm-GeneralDK-Conj}
Let $K$ be a field of characteristic zero and $z=(z_1, z_2, \dots, z_n)$ $n$ 
commutative free variables. Let $\cA$ be a commutative $K$-algebra and $V$ a $K$-subspace of $\cA$. Set $\widetilde V$ to be the $K$-subspace of the Laurent polynomial algebra $\cA[z^{-1}, z]$ consisting of the Laurent polynomials with the constant term in $V$. Then $\widetilde V$ is a  MS of $\cA[z^{-1}, z]$, if (and only if) $V$ is a MS of $\cA$.
\end{conj}

Next, we use Theorem \ref{ImDesp} to show two special cases of 
the LFED conjecture  \ref{LFED-Conj} for $\cA_n$. 

\begin{propo}\label{A1=pmI}
In terms of notation as in Lemma \ref{Lma(-2)}, assume $A_1=\pm \I_k$. 
Given $f\in \cA_n$ such that $f^m\in \im(\I-\phi)$ for all $m\ge 1$ we have 
$0\not \in \poly_z f$, 
where $\poly_z\, f$ is the polytope of $f(x)(=f(y, z))$ as 
a Laurent polynomial in $z$ over $\cB_1$.  
In particular, $\im(\I-\phi)$ is a MS of $\cA_n$.  
\end{propo}

\pf Note that the second statement follows from the first one, by 
a similar argument as in the proof of Lemma \ref{W0-Lma}, $2)$. 

To show the first statement, we assume otherwise and fix a nonzero $f(x) \in \cA_n$ 
such that $0\in \poly_z f(y, z)$ but $f^m\in \im(\I-\phi)$ for all $m\ge 1$. We drive a contradiction as follows. 

First, for all $a(y)=\sum_{\alpha\in E_1} c_\alpha y^\alpha \in \cB_1$  
the following statements can be readily verified:
{\it 
\begin{enumerate}
  \item[$i)$] if $A_1=\I_k$, then $a(y)\in \im(\I_{\cB_1}-\phi_1)$, if and only if  
$c_{\alpha}=0$ for all $\alpha\in E_1$ such that $q^\alpha=1$; 
    \item[$ii)$] if $A_1=-\I_k$, then $\phi^2=\phi$, and 
    by \cite[Proposition 5.2]{Open-LFNED} we have 
    $\im(\I_{\cB_1}-\phi_1)=\Ker(\I_{\cB_1}+\phi_1)$. Then it is easy to see that  
$a(y)\in \im(\I_{\cB_1}-\phi_1)$, if and only if $c_{-\alpha}=-q^\alpha c_{\alpha}$ 
for all $\alpha\in E_1$.
\end{enumerate}
}
Then by Theorem \ref{ImDesp} and the properties $i)$ and $ii)$ above we have
\begin{enumerate}
\item[$iii)$] {\it if $A_1=\pm I$, then 
for a Laurent polynomial 
$h(x)\in \cA_n$ we have that $h^m\in \im(\I-\phi)$ for all $m\ge 1$, 
if and only if $q_i$ $(1\le i\le n)$ and the coefficients of $h(x)$
satisfy a certain sequence of polynomial equations 
over $\bZ$. }
\end{enumerate}
 
Let $K$ be the field of fractions of the subring of $R$ 
generated by $q_i$ $(1\le i\le n)$ and the coefficients of $f(x)$. 
Then by Property $iii)$ above and a reduction similar 
as the one in Section $4$ in \cite{FPYZ} we may assume that $K$  
is a subfield of the algebraic closure $\bar \bQ$ of $\bQ$. 
Therefore, for each prime $p\ge 2$ the $p$-valuation $\nu_p$ of $\bQ$ 
can be extended to $K$, which we will still denote by $\nu_p$.

Applying Theorems \ref{DvK-Thm} and \ref{ImDesp}  we may assume  
\begin{align}\label{A1=pmI-peq1}
f(x)&=a_0(y)+\sum_{0\ne \beta\in E_2} 
 a_\beta(y) z^\beta, \\
a_0(y)&=\sum_{\alpha\in E_1} c_\alpha y^\alpha
\end{align} 
such that $a_0(y)\ne 0$ and $a_0(y) \in \im(\I_{\cB_1}-\phi_1)$.  
  
Note that by Properties $i)$ and $ii)$ $a_0(y)$ can not 
be a constant in $K$, so we can fix a nonzero extremal 
point $\lambda\in E_1$ of $\poly_y a_0(y)$. 
Let $d$ be an {\bf even} positive integer such that for all $m\ge 1$ 
\begin{align}
q^{dm\lambda}=1, \label{A1=pmI-peq2}
\end{align}
\begin{align}\label{A1=pmI-peq3-0}
c_{-\lambda}^{md}&\ne -q^{md\alpha} c_\lambda^{md}.
\end{align}

Note that such an integer $d$ does exist. This is because that, 
first,  
$q_i$ $(1\le i\le k)$  by Lemma \ref{Lma(-2)}, $1)$ 
are all roots of unity and, second, $c_\lambda\ne 0$ and, third,  
by Property $ii)$ above $c_{-\lambda}=-q^{\lambda} c_\lambda$. 

Therefore, by Properties $i)$ 
and $ii)$ above we have for all $m\ge 1$ 
\begin{align}\label{A1=pmI-peq3}
a_0(y)^{md}&\not \in \im(\I_{\cB_1}-\phi_1).
\end{align}

%
%

By Dirichlet's prime number theorem there 
exist infinitely many $m\ge 1$ such that $p\!:=md-1$ is prime.
Furthermore, it is well-known in Algebraic Number Theory 
(e.g., see \cite[Theorem 4.1.7]{W}) that for all but finitely such primes $p$,  
the values of $\nu_p$ at $q_i$ $(1\le i\le k)$ and 
the nonzero coefficients of $f(x)$ are equal to $0$.  
For all primes $p$ with these properties we consider 
\begin{align}\label{A1=pmI-peq4}
f^p(x)\equiv a_0^p(y)+\sum_{0\ne \beta\in E_2} 
 a_\beta^p(y) z^{p\beta} \mod S_p [x^{-1}, x],
\end{align}
where $S_p$ is the subring of $K$ formed by the elements 
$u\in K$ such that $\nu_p(u)\ge 1$.  

Choosing $p$ large enough we assume $p\beta+\gamma \ne 0$ for all 
$0\ne \beta, \gamma \in  \text{Poly}_z\, f(y, z)$. 
Then by Eq.\,(\ref{A1=pmI-peq4}) 
the constant term of 
$f^{p+1}$ (viewed as a Laurent polynomial in 
$\cB_1[z^{-1},z]$) is equal to $a_0^{p+1}$ modular 
$S_p[y^{-1}, y]$. 

By the choice of $p$ we have 
$p+1=md$. Since $f^{p+1}(x)=f^{md}(x)\in \im(\I-\phi)$,
by Eq.\,(\ref{A1=pmI-peq3}) and Theorem \ref{ImDesp} 
there exists a nonzero $b_{p+1}(x)\in  S_p[x, x^{-1}]$ 
such that 
\begin{align}
a_0^{p+1}(x)+ b_{p+1}(x)&\in \im(\I_{\cB_1}-\phi_1). \label{A1=pmI-peq5}
\end{align}

Since $\lambda$ is an extremal point 
of $\poly_y a_0(y)$, we have that 
$(p+1)\lambda$ is an extremal point of 
$\poly_y a_0^{p+1}(y)$, and the coefficient of 
$x^{(p+1)\lambda}$ in $a_0^{p+1}(y)$ is equal to $c_\lambda^{p+1}$.  
Then by Eqs.\,(\ref{A1=pmI-peq2}), (\ref{A1=pmI-peq3-0}), (\ref{A1=pmI-peq5}),  
and Property $i)$, $ii)$ above we have $c_\lambda^{p+1}+u_{p+1}=0$, where 
$u_{p+1}$ is the coefficient of $x^{(p+1)\lambda}$ in $b_{p+1}(x)$. 
Hence $0\ne c_\lambda^{p+1}=-u_{p+1}\in S_p$, which implies  
$\nu_p (c_\lambda^{p+1})\ge 1$. 
Therefore we also have  
$\nu_p(c_\lambda)\ge 1$. But this contradicts to the choice 
of $p$ with $\nu_p(c_\lambda)=0$. 
\epfv

It is worthy to point out explicitly the following two special cases of the proposition above. 

\begin{corol}\label{ScalarCase}
Assume $A=\pm I$, i.e., $\phi$ is the $R$-algebra automorphism of $\cA$ that maps $x_i$ to $q_ix_i$ for all $1\le i\le n$, or $\phi$ maps $x_i$ to $q_ix_i^{-1}$ for all $1\le i\le n$. Then $\im(\I-\phi)$ is a MS of $\cA$.
\end{corol}

Next, we show the LFEN conjecture (and also the LNED conjecture)
for the Laurent polynomial algebras 
over a field of characteristic zero in one or 
two commutative free variables.  

\begin{theo}\label{One-Two-VariableCase}
Let $K$ be a field of char zero and $\cA$ the Laurent polynomial algebra over $K$ in one or two commutative variables. Then both the LFED Conjecture \ref{LFED-Conj} and LNED Conjecture \ref{LNED-Conj} holds for $\cA$.
\end{theo}

\pf By Theorem \ref{LNED-Case} we  only need to show 
the $K$-$\cE$-derivation case for the LFED 
Conjecture \ref{LNED-Conj}. Furthermore, 
the one variable case follows from Corollary \ref{ScalarCase} 
above. So we may assume that $\cA=\cA_2$ with $R=K$, 
i.e., $\cA$ is the Laurent polynomial algebra over $K$  
in two variables.

Let $\phi$ be a LF $K$-algebra endomorphism of $\cA_2$. 
We use the same notations as before, in particular, 
let $A$ be the exponent matrix of $\phi$. 
If the rank of $A$ is less than $2$, then by Lemmas \ref{Lma5.2.1}, \ref{FuRankRed} and Corollary \ref{ScalarCase} it is easy to see that 
$\im(\I-\phi)$ also is a MS of $\cA_2$. 

So we assume that $A$ is of rank $2$.  
Then the case $k=2$, i.e., $E_2=\bZ^2$, follows from Lemma \ref{URootCase}; 
the case $k=1$ follows from Proposition \ref{A1=pmI}; and the case $k=0$ follows from Lemma \ref{W0-Lma}. 
\epfv
 
We end this section with the following corollary of the proof of 
Theorem \ref{ImDesp}, which can be proved by going through the proof of Theorem \ref{ImDesp} with some slight modifications for the polynomial algebra $R[x]$ instead of 
$\cA_n=R[x^{-1}, x]$. 

\begin{corol}\label{PolyAlgCase}
Let $\phi$ be as in Theorem \ref{ImDesp}, and assume further that 
all entries of the exponent matrix $A$ of $\phi$   
lie in $\bN$.  Then the following statements hold: 
\begin{enumerate}
  \item[$1)$] $\phi$ preserves the polynomial algebra $R[x]$ 
  and its restriction $\psi\!:=\phi\,|_{R[x]}$ 
      is a locally finite $R$-derivation of $R[x]$.
  \item[$2)$] $\im(\I_{R[x]}-\psi)$ consists of all $f\in R[x]$ of the form 
 \begin{align} 
f(x)=a_0(y)+\sum_{0\ne \beta\in E_2\cap \bN^n} 
 a_\beta(y) z^\beta
\end{align}
with $a_\beta(y)\in R[y]$ $(\beta\in E_2\cap\bN^n)$ and 
$a_0(y)\in \im(\I_{R[y]}-\phi_{R[y]})$.  
\end{enumerate} 

In particular, $\im(\I_{R[x]}-\psi)$ is a MS of $R[x]$, i.e., 
the LFED conjecture \ref{LFED-Conj} holds for the $R$-$\cE$-derivation 
$\psi$ of $R[x]$. 
\end{corol}

\end{document}